\documentclass[11pt]{amsproc}
\usepackage{diagrams}
\def\thmname{THEOREME}
\def\propname{PROPOSITION}
\def\lemmaname{LEMME}
\def\corname{COROLLAIRE}
\def\exname{Exemple}
\def\remname{Remarque}

\def\d{\partial}

\def\dsum{\mathop{\displaystyle \sum }}
\def\stackunder#1#2{\mathrel{\mathop{#2}\limits_{#1}}}

\font\bb=msbm10
\def\G{\mathcal G}
\def\R{\hbox{\bb R}}

\def\rarrow{\mathop{\hbox to 4cm{\rightarrowfill}}}
\def\larrow{\mathop{\hbox to 4cm{\leftarrowfill}}}
\newcommand\object[1] {\makebox(0,0)[c]{$#1$}}
\def\QED {\vrule height 1.2ex width 1.2ex}{\vskip 10pt plus100pt}

\title{\sc Une formule du type Baker-Campbell-Hausdorff pour les 
groupo{\"\i}des de Lie}
\author{\sc Birant Ramazan\thanks{Work partially supported by the Grant 
104/1999 of the Romanian Academy}}

\begin{document}

\newtheorem{thm}{\thmname}[section]
\newtheorem{prop}[thm]{\propname}
\newtheorem{lem}[thm]{\lemmaname}
\newtheorem{defn}[thm]{\defnname}
\newtheorem{cor}[thm]{\corname}
\newtheorem{ex}[thm]{\exname}
\newtheorem{rem}[thm]{\remname}

\begin{abstract}
On d\'emontre dans le contexte d'un groupo{\"\i}de de Lie $G$ un analogue 
de la formule de Baker-Campbell-Hausdorff. Comme application on calcule les 
fonctions de structure de l'alg\'ebro{\"\i}de de Lie associ\'e \`a $G$.
\end{abstract}

\maketitle

\section{INTRODUCTION}
Cet article est consacr\'e a l' \'etude locale d'un groupo{\"\i}de de Lie $G$
dans des cartes convenablement choisies. En particulier on obtient 
dans le th\'eor\`eme \ref{Taylorgroupoide} le developpement de la 
multiplication ce qui constitue un analogue de la formule de 
Baker-Campbell-Hausdorff du cas des groupes de Lie (cf. \cite{kirillov} 
par exemple). Ce r\'esultat a \'et\'e pr\'esent\'e bri\`evement dans
\cite{ramazan1}, \cite{ramazan2} ou il est utilis\'e pour demontrer que le
groupo{\"\i}de tangent $\tilde{G}$ associ\'e \`a un groupo{\"\i]de de Lie $G$
est lui aussi un groupo{\"\i}de de Lie. Il intervient aussi de mani\`ere 
essentielle dans le calcul du commutateur dans l'alg\`ebre de convolution du 
groupo{\"\i}de tangent $\tilde{G}$, ce qui permet de quantifier la structure 
de Poisson canonique du dual de l'alg\'ebro{\"\i}de de Lie associ\'e \`a $G$.

L'article est structur\'e comme il suit. Apr\`es avoir fix\'es la terminologie
et les notations, on rappelle pour le benefice du lecteur les diff\'erentes 
constructions de l'alg\'ebro{\"\i}de de Lie ${\mathcal G}$ associ\'e \`a un 
groupo{\"\i}de de Lie $G$. Dans la suite on explicite la structure locale de 
$G$ dans une carte et on montre comment on associe \`a une carte de $G$ 
choisie convenablement, une carte de son alg\'ebro{\"\i}de de Lie ${\mathcal G}$. 
On \'ecrit la multiplication et l'inversion de $G$ dans ces cartes et on les 
developpe en s\`eries de Taylor pour obtenir l'analogue de la formule de 
Baker-Campbell-Hausdorff. Enfin, comme application, on calcule les fonctions 
de structure de ${\mathcal G}$. 

Rappelons bri\`evement les principaux faits sur les groupo{\"\i}des de Lie
et les alg\'ebro{\"\i}des de Lie associ\'es. Pour une pr\'esentation 
detaill\'ee de la th\'eorie des groupo{\"\i}des de Lie on pourra se reporter 
aux ouvrages de A.~Weinstein, P.~Dazord et A.~Coste \cite{coste}, o\`u 
K.~Mackenzie \cite{mackenzie}. Les notations et les d\'efinitions de la 
th\'eorie des groupo{\"\i}des seront celles donn\'ees dans \cite{renault1} par 
J.~Renault. Par d\'efinition un groupo{\"\i}de est un ensemble $G$ muni d'un 
produit $G\times G\supset G^{(2)}\ni (g_1,g_2)\longmapsto g_1g_2\in G$ 
d\'efini sur un sous-ensemble $G^{(2)}$ de $G\times G$, et une application 
inverse $G\ni g\longmapsto g^{-1}\in G$ v\'erifiant:

\smallskip

1. $(g^{-1})^{-1}=g$

\medskip

2. Si $(g_1,g_2),(g_2,g_3)\in G^{(2)}$ alors $(g_1g_2,g_3), (g_1,g_2g_3)\in 
G^{(2)}$ et $(g_1g_2)g_3=g_1(g_2g_3)$

\medskip

3. $(g^{-1},g)\in G^{(2)}$. Si $(g_1,g_2)\in G^{(2)}$ alors $g_1^{-1}(g_1g_2)=
g_2$

\medskip

4. $(g,g^{-1})\in G^{(2)}$. Si $(g_2,g_1)\in G^{(2)}$ alors $(g_2g_1)g_1^{-1}=
g_2$

\medskip

\noindent $G^{(2)}$ s'appelle l'ensemble des paires composables et pour 
$g\in G$ on appelle $s(g)=g^{-1}g$ le domaine de $g$ et $r(g)=gg^{-1}$ l'image 
de $g$. La composition $g_1g_2$ est bien d\'efinie si et seulement si 
$r(g_2)=s(g_1)$. L'ensemble $s(G)=r(G)$, not\'e $G^{(0)}$, sera identifi\'e 
\`a une partie de $G$ et appel\'e espace des unit\'es. Pour $x\in G^{(0)}$, 
on notera $G^x=r^{-1}(x)$, $G_x=s^{-1}(x)$.

Un groupo{\"\i}de de Lie est un groupo{\"\i}de $G$ qui a une structure 
de vari\'et\'e diff\'erentiable compatible avec la structure alg\'ebrique :
\begin{enumerate}
\item $G^{(0)}$ est une sous-vari\'et\'e de $G$
\item $r,s:G\rightarrow G^{(0)}$ sont des submersions
\item la multiplication : $G^{(2)}\rightarrow G$ est diff\'erentiable
\end{enumerate}

Comme consequences de la d\'efinition il faut noter que l'application 
$i:G\longrightarrow G$, $i(\gamma )=\gamma ^{-1}$ est un diff\'eomorphisme 
(voir \cite{mackenzie}, p.85), et aussi le fait qu'en notant 
$m=\mbox{dim}\, G$ et $n=\mbox{dim}\, G^{(0)}$, pour tout $x\in G^{(0)}$, 
$G^x$ et $G_x$ sont des sous-vari\'et\'es de $G$ de dimension $m-n$.

On rappelle maintenant les differentes constructions de l'alg\'ebro{\"\i}de de 
Lie associ\'e \`a un groupo{\"\i}de de Lie $G$ de base $G^{(0)}$. Les 
alg\'ebro{\"\i}des de Lie ont \'et\'e introduits par J.~Pradines 
\cite{pradines2}, et g\'en\'eralisent la notion d'alg\`ebre de Lie dans le 
cadre de la th\'eorie des groupo{\"\i}des de Lie. 

Pour fixer les notations, pour toute application diff\'erentiable $f$ entre 
les vari\'et\'es $M$ et $N$, $Tf$ designe l'application tangente et $T_xf$ 
l'application tangente en $x$ entre les espaces tangents $T_xM$ et $T_{f(x)}N$.
Aussi pour $E$ fibr\'e vectoriel de classe $C^{\infty}$ sur la vari\'et\'e $M$ 
on notera par $C^{\infty}(M,E)$ l'ensemble des sections de classe $C^{\infty}$ 
de $E$ sur $M$.

Par d\'efinition un alg\'ebro{\"\i}de de Lie sur une vari\'et\'e $M$ est un 
triplet constitu\'e d'un fibr\'e vectoriel $E$ de base $M$ et classe $C^{
\infty }$, une structure de {\R}-alg\`ebre de Lie sur l'espace des sections
$C^{\infty}(M,E)$, dont on note $\left[ \cdot ,\cdot \right] $ le crochet et 
un morphisme $\rho :E\rightarrow TM$ de fibr\'es vectoriels $C^{\infty}$, 
appel\'e ancre tels que:

(i) L'application induite entre les espaces des sections $\overline{\rho }: 
C^{\infty} 
(M,E) \rightarrow C^{\infty} (M,TM)$, $\overline{\rho}(\xi )(x)=\rho (\xi 
(x))$, $\xi \in C^{\infty}(M,E)$,\ $x\in M$,  est un morphisme d'alg\`ebres 
de Lie : $\left[ \overline \rho (\xi ),\overline \rho (\eta )\right] =
\overline \rho \left( \left[ \xi , \eta \right] \right)$

(ii) Pour toute fonction $f\in C^{\infty }(M)$ et pour tout couple 
$(\xi ,\eta )$ de sections $C^{\infty }$ de $E$ ,

\centerline{$\left[ \xi,f\eta \right] = f\left[ \xi ,\eta \right] +
\overline{\rho }(\xi )(f)\eta $}

\medskip

On fixe $\left\{e_1,e_2,...,e_p\right\}$ un rep\`ere local sur $U\subset M$ 
pour $E$ et $(q_1,...,q_n,
\lambda _1,...,\lambda _p)$ coordonn\'ees locales de $E$ avec les $q_i$ 
coordonn\'ees locales pour la base $M$ et les $\lambda_j$ coordonn\'ees dans 
les fibres associ\'ees au rep\`ere $\left\{ e_1,e_2,...,e_p\right\}$.
Alors, localement, le fait que $E$ est un alg\'ebro{\"\i}de de Lie implique 
l'existence des fonctions de structure $c_{ijk}, a_{ij}\in C^{\infty}(U)$ 
telles que $\left[ e_i,e_j\right] =\stackunder{k}{\dsum}c_{ijk}e_k$ et
$\overline{\rho}(e_i)=\stackunder{j}{\dsum}a_{ij}\dfrac{\d}{\d q_j}$.

Remarquons tout d'abord que pour tout $g\in G$, $R_g:G_{r(g)}\rightarrow 
G_{s(g)}$, $R_g(h)=hg$ et $L_g:G^{s(g)}\rightarrow G^{r(g)}$, $L_g(h)=gh$ sont 
des diff\'eomorphismes et que pour tout $g\in G$ on a $T_gG^{r(g)}=KerT_gr$ et 
$T_gG_{s(g)}=KerT_gs$. Cela permet de d\'efinir les champs invariants \`a 
gauche sur $G$ par $L(G)=\left\{ \xi \in C^{\infty }(G,TG)\mid \xi \in KerTr, 
TL_g\circ \xi =\xi \circ L_g\right\}$ et les champs invariants \`a droite par 
$R(G)=\left\{ \xi \in C^{\infty }(G,TG)\mid \xi \in KerTs, TR_g\circ \xi =\xi 
\circ R_g\right\}$. Il est facile a voir que $L(G)$ et $R(G)$ sont des 
alg\`ebres de Lie.

Le fait que $G^{(0)}$ est une sous-vari\'et\'e de $G$ permet de considerer 
$T_xG^{(0)}$ sous-espace de $T_xG$, pour tout $x\in G^{(0)}$. 
On note ${\mathcal L}$, ${\mathcal R}$, respectivement ${\mathcal N}$, les fibr\'es 
vectoriels sur $G^{(0)}$ dont les fibres au dessus de $u\in G^{(0)}$ sont 
${\mathcal L}_u=KerT_ur$, ${\mathcal R}_u=Ker T_us$, respectivement
${\mathcal N}_u=T_uG/T_uG^{(0)}$.

\begin{lem}\label{217}
On a les isomorphismes des espaces vectoriels $L(G)\simeq C^{\infty}(G^{(0)},
{\mathcal L})$ et $R(G)\simeq C^{\infty}(G^{(0)},{\mathcal R})$.
\end{lem}

{\it Preuve.} L'application $\Phi :L(G)\rightarrow C^{\infty }(G^{(0)},
{\mathcal L}),\Phi (\xi )=
\xi _{\mid _{G^{(0)}}}$ est bien d\'efinie.
Montrons que $\Phi $ est injective. Soit $\xi _{\mid _{G^{(0)}}}=\eta _{\mid 
_{G^{(0)}}}$. Alors pour $g\in G$ on a $\xi (g)=\xi (gu)=T_uL_g\xi (u)=
T_uL_g\eta (u)=\eta (gu)=\eta (g)$, o\`u $u=s(g)$.

Montrons que $\Phi $ est surjective. Soit $\eta \in C^{\infty }(G^{(0)},{\mathcal 
L})$. On d\'efinit $\xi :G\longrightarrow TG,$ $\xi(g)=T_uL_g\eta (u)$, pour 
$u=s(g)$. On voit que $L_u=Id_{G^u}$, d'o\`u $\xi (u)=\eta (u)$ pour $u\in 
G^{(0)}$, donc 
$\eta =\xi \mid _{G^{(0)}}$. Il reste \`a montrer que $\xi \in L(G)$. On a 
$T_{g^{\prime }}L_g\xi (g^{\prime })=T_{g^{\prime }}L_g\left[ T_uL_{g'}\xi (u)
\right] =T_uL_{gg^{\prime }}\xi (u)=\xi (gg^{\prime })$, 
o\`u $u=s(g^{\prime })$. On a utilis\'e le fait que $T_{g^{\prime }}L_g\circ 
T_uL_{g^{\prime }}=T_u(L_g\circ L_{g^{\prime }})=T_uL_{gg^{\prime }}$. 
Comme $\Phi $ est \'evidemment lin\'eaire on a d\'emontr\'e le premier 
isomorphisme. La d\'emonstration du second isomorphisme est analogue.
\QED

\medskip

Les fibr\'es ${\mathcal L}$, avec le crochet donn\'e par $\left[ \xi ,\eta 
\right]=\Phi (\left[ \Phi ^{-1}(\xi ),\Phi ^{-1}(\eta )\right] )$ et l'ancre 
$\rho :{\mathcal L} \rightarrow TG^{(0)}$, $\rho _u=T_us$, respectivement 
${\mathcal R}$, avec le crochet d\'efini de mani\`ere analogue \`a ${\mathcal 
L}$ et l'ancre $\mu :{\mathcal R}\rightarrow TG^{(0)}$, $\mu _u=T_ur$, sont 
deux alg\'ebro{\"\i}des de Lie antiisomorphes par l'application tangente $Ti$ 
de l'inversion $i$ de $G$.

L'application $I_x-T_xr:T_xG\rightarrow Ker T_xr$, $X\mapsto X-T_xrX$, est
bien d\'efinie et surjective. Son noyau est $T_xG^{(0)}$ et par factorisation 
on obtient l'isomorphisme d'espaces vectoriels ${\mathcal N}_x\simeq Ker T_xr$.
De la m\^eme mani\`ere, en consid\`erant $I_x-T_xs:T_xG\rightarrow Ker T_xs$, 
$X\mapsto X-T_xsX$ on d\'emontre l'isomorphisme ${\mathcal N}_x\simeq Ker 
T_xs$.

Ces deux isomorphismes d\'efinissent sur ${\mathcal N}$ deux structures 
d'alg\'ebro{\"\i}de de Lie antiisomorphes. Le crochet de Lie sur 
$C^{\infty}(G^{(0)},{\mathcal N})$ est d\'efini en utilisant
l'isomorphisme $L(G)\ni \xi \mapsto \left[ \xi _{G^{(0)}}\right] \in 
C^{\infty}(G^{(0)},{\mathcal N})$, o\`u $\left[ X\right]$ est l'image de 
$X\in T_xG$ dans $T_xG/T_xG^{(0)}$. L'ancre sur ${\mathcal N}$ est 
$\rho :{\mathcal N}\rightarrow TG^{(0)}$, $\rho _x=T_xr-T_xs$. 

Dans la suite on appellera alg\'ebro{\"\i}de de Lie du groupo{\"\i}de 
de Lie $G$ le fibr\'e ${\mathcal L}$ avec la structure d'alg\'ebro{\"\i}de 
d\'efinie pr\'ec\'edemment et pour mettre en \'evidence qu'il est 
l'alg\'ebro{\"\i}de associ\'e au groupo{\"\i}de $G$ il sera 
not\'e ${\mathcal G}$.

\section{La structure locale d'un groupo{\"\i}de de Lie}
\subsection{Les cartes}

Soit $G$ un groupo{\"\i}de de Lie et ${\mathcal G}$ son alg\'ebro{\"\i}de de Lie. 
On va expliciter dans cette section la structure de $G$ dans une carte 
convenablement choisie au voisinage d'un point $x_0\in G^{(0)}\subset G$.
Des cartes de ce genre ont \'et\'e utilis\'ees aussi dans \cite{nistor}.

Comme $r$ est une submersion au point $x_0$ appartenant \`a la 
sous-vari\'et\'e $G^{(0)}$ de $G$, il existe $U$ voisinage ouvert de 0 dans 
$\R^n$, $V$ voisinage ouvert de 0 dans $\R^m$ et les cartes $\psi :U\times 
V\rightarrow G$, $\varphi :U\rightarrow G^{(0)}$ v\'erifiant:

\begin{equation}\label{psibon}
\begin{array}{l}
1.\ \psi (0,0)=x_0\\
2.\ r(\psi (u,v))=\varphi (u)\\
3.\ \psi (U\times \lbrace 0\rbrace )=\psi (U\times V)\bigcap G^{(0)}
\end{array}
\end{equation}

\noindent La deuxi\`eme condition revient au diagramme commutatif :

\bigskip

\centerline{
\setlength \unitlength {0.25em}
\begin{picture}(150,30)(-75,-15)
\put (-40, -6){\vector(0,1){12}}
\put (40, -6){\vector(0,1){12}}
\put (-25, 10){\vector(1,0){52}}
\put (-30, -10){\vector(1,0){60}}
\put(-40, 10){\object{G\supset \psi (U\times V)}}
\put(-40, -10){\object{U\times V}}
\put(40, 10){\object{\varphi (U)\subset G^{(0)}}}
\put(40, -10){\object{U}}
\put(45, 0){\object{\varphi}}
\put(-45, 0){\object{\psi }}
\put(0, 13){\object{r}}
\put(0, -7){\object{pr_1}}
\end{picture}}

\bigskip

Des deux derni\`eres conditions on d\'eduit $\varphi (u)=\psi (u,0)$, et en 
cons\'equence on pourra exprimer la structure de $G$ en utilisant seulement 
la carte $\psi$. Toutefois pour la simplicit\'e des notations on gardera 
$\varphi =\psi(\cdot , 0)$.

A la carte $\psi$ de $G$ s'associe canoniquement une carte de 
l'alg\'ebro{\"\i}de de Lie ${\mathcal G}$. Plus pr\'ecis\'ement on a :

\begin{lem}\label{cartealg}
L'application $\theta :U\times \R^m\rightarrow {\mathcal G}$, $\theta (u,v)=
(\varphi (u),\dfrac{\d \psi}{\d v}(u,0)v)$ est une carte de ${\mathcal G}$ au 
voisinage de la fibre ${\mathcal G}_{x_0}$ et la famille $\lbrace e_1,e_2,..,e_m
\rbrace$ d\'efinie par $e_i(\varphi (u))=\theta (u,f_i)$, $i=\overline{1,m}$, 
o\`u $\lbrace f_1,f_2,...,f_m\rbrace $ est la base canonique de $\R^m$, est
un rep\`ere mobile de ${\mathcal G}$ sur $\varphi (U)$.
\end{lem}

{\it Preuve.} Pour tout $u\in U$ on a 

\centerline{$G^{\varphi (u)}\bigcap \psi (U\times V)=\lbrace \gamma \in \psi 
(U\times V)|r(\gamma )=\varphi (u)\rbrace =\psi (\lbrace u\rbrace \times V)$} 

\noindent L'application $\psi (u,\cdot ):V\rightarrow G^{\varphi (u)}$ est 
alors une carte de la sous-vari\'et\'e $G^{\varphi (u)}$, qui associe $0\in V$ 
\`a $\varphi (u)$. On peut identifier $\lbrace u\rbrace \times \R^m $ avec $T_{
\varphi (u)}G^{\varphi (u)}$ par l'isomorphisme $\dfrac{\d \psi}{\d v}(u,0):T_{
(u,0)}(\lbrace u\rbrace \times V)=\lbrace u\rbrace \times \R^m \rightarrow 
T_{\varphi (u)}G^{\varphi (u)}$.
L'image de $\theta$ est le voisinage $\theta (U\times \R^m)={\mathcal G}_{\varphi 
(U)}$ de la fibre ${\mathcal G}_{x_0}$. \QED

\bigskip

{\it Remarque.}
Dans un groupe de Lie il existe un voisinage de l'unit\'e diff\'eomorphe avec 
un voisinage de l'\'el\'ement nul de l'alg\`ebre de Lie associ\'e. Le lemme 
pr\'ec\'edent permet de donner la g\'en\'eralisation suivante pour les 
groupo{\"\i}des de Lie~:

Pour tout $x_0\in G^{(0)}$ il existe un voisinage de $x_0$ dans $G$ qui est 
diff\'eomorphe avec un voisinage de $(x_0,0)$ dans $\G$. 

\setlength \unitlength {0.25em}
\begin{picture}(130,30)(-65,-15)
\put (-21, 10){\vector(1,0){42}}
\put(-35, 10){\object{{\G}\supset \theta (U\times V)}}
\put(35, 10){\object{\psi (U\times V)\subset G}}
\put(0, -15){\object{U\times V}}
\put(-3,-11){\vector(-3,2){27}}
\put(3,-11){\vector(3,2){28}}
\put(22, -5){\object{\psi}}
\put(-22, -5){\object{\theta}}
\put(0, 13){\object{\alpha}}
\put(0, 7){\object{\widetilde{}}}
\end{picture}

En effet avec les notations pr\'ecedentes, $\alpha =\psi \circ \theta ^{-1}$ 
est un diff\'eomorphisme entre le voisinage $\theta (U\times V)$ de 
$(x_0,0)\in \G$ et le voisinage $\psi (U\times V)$ de $x_0$ dans $G$. On 
remarque de plus que pour tout $x\in \varphi (U)$, $\alpha (\G _x)\subset G^x$.

Ce r\'esultat n'est qu'un cas particulier de la proposition suivante qui
est bas\'ee sur l'existence d'une application exponentielle pour tout
groupo{\"\i}de de Lie. Cette application exponentielle introduite par Pradines 
dans \cite{pradines3} g\'en\'eralise \`a la fois l'exponentielle d'un 
groupe de Lie et l'exponentielle d'une vari\'et\'e munie d'une connexion.

\begin{prop}\label{diffeoG}
Soit $G$ un groupo{\"\i}de de Lie et $\G$ son alg\'ebro{\"\i}de de Lie. Il
existe alors un voisinage $V$ de $G^{(0)}$ vu comme la section nulle 
$\left\{ (x,0)|x\in G^{(0)}\right\}$ dans $\G$, un voisinage $W$ de 
$G^{(0)}$ dans $G$ et un diff\'eomorphisme $\alpha :V\rightarrow W$ tel que
$\alpha (\G _x \bigcap V)=G^x\bigcap W$ et $\alpha ^{\prime}_x (0)$ est 
l'identit\'e de $\G _x$, o\`u $\alpha _x$ est la restriction de 
$\alpha$ sur $\G _x\bigcap V$.
\end{prop}

\noindent L'id\'ee de la d\'emonstration est la suivante. Soit $\nabla$ 
une connexion sur l'alg\'ebro{\"\i}de de Lie $\G$. On associe \`a $\nabla$ 
une connexion invariante \`a gauche sur $G$, dont la restriction \`a $G^x$ 
est une connexion lin\'eaire $\nabla _x$. On peut alors d\'efinir fibre par 
fibre une application exponentielle, et prendre comme $\alpha$ cette 
exponentielle. Pour les d\'etails voir \cite{landsman5} ou \cite{nistor}.

\subsection{La multiplication et l'inversion}

Pour exprimer le produit et l'inversion de $G$ dans la carte $\psi$ on 
a besoin de la forme de l'application source $s$ dans cette carte, forme qui 
est explicit\'ee dans le lemme suivant.

\begin{lem}
Il existe une submersion $\sigma :U\times V\rightarrow U$ telle que 
$s(\psi (u,v))=\varphi (\sigma (u,v))$. De plus $\sigma (u,0)=u$.
\end{lem}

{\it Preuve.}
En r\'eduisant eventuellement $V$, on peut supposer que $s(\psi (u,v))\in 
\varphi (U)$, pour $(u,v)\in U\times V$. Il existe alors un \'el\'ement 
$\sigma (u,v)\in U$ tel que $s(\psi (u,v))=\varphi (\sigma (u,v))$. Evidemment 
$\sigma =\varphi ^{-1}\circ s\circ \psi$ est une submersion, comme expression 
dans les cartes de la submersion $s$. Enfin $\varphi (\sigma (u,0))=
s(\psi (u,0))=\psi (u,0)=\varphi (u)$, donc $\sigma (u,0)=u$. \QED

On peut maintenant donner les d\'eveloppements dans la carte $\psi$ de la 
multiplication et de l'inversion de $G$. Le r\'esultat suivant repr\'esente 
l'analogue de la formule de Baker-Campbell-Hausdorff pour les groupo{\"\i}des 
de Lie.

\begin{prop}\label{Taylorgroupoide}
(i) Pour $u,u_1\in U$ et $v,w\in V$ on a 
$(\psi (u,v),\psi (u_1, w))\in G^{(2)}$ si et seulement si 
$u_1=\sigma (u,v)$. Dans ce cas le produit est donn\'e par
$\psi (u,v)\psi (\sigma (u,v),w) =\psi (u,p(u,v,w))$
o\`u $p:U\times V\times V\rightarrow V$ est une application
diff\'erentiable qui a un d\'eveloppement de la forme 
$p(u,v,w)=v+w+B(u,v,w)+O_3(u,v,w)$
avec $B$ bilin\'eaire en $(v,w)$ et $O_3(u,v,w)$ de l'ordre de 
$\| (v,w)\|^3$.

(ii) Soit $(u,v)\in U\times V$ tel que $\psi (u,v)^{-1}\in \psi (U\times 
V)$. Alors $\psi (u,v)^{-1}=\psi (\sigma (u,v),w)$, o\`u $w$ v\'erifie 
$p(u,v,w)=0$. De plus on a le d\'eveloppement $w=-v+B(u,v,v)+O_3(u,v)$, avec
$O_3(u,v)$ de degr\'e d'homog\'en\'eit\'e superieur \`a 3 en $v$.
\end{prop}

{\it Preuve.}
{(i)} Soit $g=\psi (u,v)$ et $h=\psi (u_1,w)$. On a $s(g)=\varphi (\sigma 
(u,v))$ et $r(h)=\varphi (u_1)$ ce qui montre que $(g,h)\in G^{(2)}$ si et 
seulement si $u_1=\sigma (u,v)$. De plus $r(gh)=r(g)=\varphi (u)$ assure 
l'existence d'un unique $p(u,v,w)\in V$ tel que

\centerline{$\psi (u,v)\psi (\sigma (u,v),w)=\psi (u,p(u,v,w))$}

\noindent On d\'efinit ainsi l'application $p:U\times V\times V\rightarrow V$, 
qui v\'erifie en particulier $p(u,0,w)=w$ et $p(u,v,0)=v$. En effet 

$\psi (u,p(u,0,w))=\psi (u,0)\psi (\sigma (u,0),w)=\varphi (u)\psi (u,w)=\psi 
(u,w)$

$\psi (u,p(u,v,0))=\psi (u,v)\psi (\sigma (u,v),0))=
\psi (u,v)s(\psi (u,v))=\psi (u,v)$.

\noindent Il s'ensuit que $\dfrac{\d p}{\d v}(u,v,0)=I$, $\dfrac{\d p}{\d w}
(u,0,w)=I$, $\dfrac{\d ^2 p}{\d v^2}(u,0,0)=0$, $\dfrac{\d ^2 p}{\d w^2}
(u,0,0)=0$, et par un d\'eveloppement de Taylor $p(u,v,w)=p(u,0,0)+
\dfrac{\d p}{\d v}(u,0,0)
v+\dfrac{\d p}{\d w}(u,0,0)w+B(u,v,w)+O_3(u,v,w)=v+w+B(u,v,w)+
O_3(u,v,w)$, o\`u $B(u,v,w)$ est pour chaque $u$ bilin\'eaire en $(v,w)$ et 
$O_3(u,v,w)$ est homog\`ene d'un degr\'e superieur \`a 3 en $v$ et $w$.

\smallskip

{(ii)} Soit $g=\psi (u,v)$. On cherche $u_1$ et $w$ tels que $g^{-1}=
\psi (u_1,w)$. D'une part $r(g^{-1})=\varphi (u_1)$ et $s(g)=\varphi (\sigma 
(u,v))$ impliquent $u_1=\sigma (u,v)$. D'autre part comme $gg^{-1}=r(g)$, on 
d\'eduit $\psi (u,v)\psi (\sigma (u,v),w)=\psi (u,0)$, donc $p(u,v,w)=0$.
Le th\'eor\`eme des fonctions implicites assure l'existence d'un $f$ 
diff\'erentiable tel que $w=f(u,v)$. On d\'eveloppe 

$f(u,v)=f(u,0)+\dfrac{\d f}{\d v}(u,0)v+...=f_1(u,v)+f_2(u,v)+...$ 

\noindent o\`u $f_k(u,v)$ est de degr\'e d'homog\'en\'eit\'e $k$ en $v$. 
On va d\'eterminer $f_1$ et $f_2$. Pour cela on utilise $p(u,v,w)=0$. On a donc

$0=p(u,v,w)=v+w+B(u,v,w)+...$

\hspace*{0,38cm}$=v+f_1(u,v)+f_2(u,v)+...+B(u,v,f_1(u,v)+f_2(u,v)+...)+...$

\noindent Le terme de degr\'e d'homog\'en\'eit\'e 1 dans le d\'eveloppement 
pr\'ec\'edant est $v+f_1
(u,v)$ donc $f_1(u,v)=-v$. On remplace dans l'\'equation pr\'ec\'edente et on 
trouve $0=f_2(u,v)+B(u,v,-v)+...$. En identifiant le terme de degr\'e 2 on 
obtient $f_2(u,v)=B(u,v,v)$. \QED

\section{Le calcul des fonctions de structure de ${\mathcal G}$}

On se propose de calculer les fonctions de structure de l'alg\'ebro{\"\i}de de 
Lie ${\mathcal G}$. On rappelle que $a_{ij}=\overline{\rho}(e_i)(q_j)$ et les 
$c_{ijk}$ sont donn\'ees par $\left[ e_i,e_j\right] =\dsum c_{ijk}e_k$, o\`u
$\rho =Ts$ est l'ancre de ${\mathcal G}$, $\left\{ e_1,...,e_m\right\}$ le 
rep\`ere mobile de ${\mathcal G}$ d\'efini dans le lemme \ref{cartealg} et 
$q_j=pr_j\circ \varphi ^{-1}$ sont les fonctions de coordonn\'ees de 
$G^{(0)}$. On notera par $B_1,..,B_m$ les coordonn\'ees 
de l'application $B:U\times V\times V\rightarrow V$
dans la base $\left\{ f_1,f_2,..,f_m\right\}$ de $\R ^m$.

\begin{prop}\label{329}
Pour tout $u\in U$, les fonctions de structure de l'alg\'ebro{\"\i}de 
${\mathcal G}$ sont donn\'ees par $a_{ij}(\varphi (u))=\dfrac{\d \sigma 
_j}{\d v_i}(u,0)$ et $c_{ijk}(\varphi (u))=B_k(u,f_i,f_j)-B_k(u,f_j,f_i)$.
\end{prop}

{\it Preuve.} {\it (i) Le calcul de } $a_{ij}$

\medskip

On obtient la forme de $a_{ij}$ par le calcul suivant 

\medskip

$a_{ij}(\varphi (u))=\overline{\rho }(e_i)(q_j)(\varphi (u))=\rho (e_i(
\varphi (u)))(q_j)=(T_{\varphi (u)}s)(e_i(\varphi (u)))(q_j)$

\medskip

\hspace*{1,8cm}$=e_i(\varphi (u))(q_j\circ s)=(\dfrac{\d \psi}{\d v}(u,0)f_i)
(pr_j\circ \varphi ^{-1}\circ s)$

\medskip

\hspace*{1,8cm}$=\dfrac{\d }{\d v_i}(pr_j\circ \varphi^{-1}
\circ s\circ \psi )(u,0)=\dfrac{\d \sigma _j}{\d v_i}(u,0)$.

\medskip

{\it (ii) Le calcul de} $c_{ijk}$

\smallskip

{\bf 1.} On explicite d'abord la forme d'une section $\xi \in C^{
\infty}(G^{(0)},{\mathcal G})$ dans les cartes de $G^{(0)}$ et ${\mathcal G}$ 
engendr\'ees par $\psi$. 

\bigskip

\centerline{
\setlength \unitlength {0.25em}
\begin{picture}(150,30)(-75,-15)
\put (-35, -6){\vector(0,1){12}}
\put (38, -6){\vector(0,1){12}}
\put (-26, 10){\vector(1,0){56}}
\put (-30, -10){\vector(1,0){60}}
\put(-40, 10){\object{G^{(0)}\supset \varphi (U)}}
\put(-40, -10){\object{\R ^n\supset U}}
\put(40, 10){\object{{\mathcal G}_{\varphi (U)}}}
\put(40, -10){\object{U\times \R ^m}}
\put(42, 0){\object{\theta}}
\put(-38, 0){\object{\varphi }}
\put(0, 13){\object{\xi }}
\put(0, -7){\object{\Xi}}
\end{picture}}

\bigskip

\noindent Cette forme est donn\'ee par $\Xi (u)=\theta ^{-1}\circ \xi \circ 
\varphi (u)=(u,\dfrac{\d \psi}{\d v}(u,0)^{-1}(\xi (\varphi (u))))$. On note 
$\xi _0:U\rightarrow \R ^m$, $\xi _0(u)=\dfrac{\d \psi}{\d v}(u,0)^{-1}(\xi 
(\varphi (u)))$ et on a $\Xi (u)=(u,\xi _0(u)$.

\medskip

{\bf 2.} Par le lemme \ref{217}, on associe \`a tout $\xi \in C^{
\infty}(G^{(0)},{\mathcal G})$ une section \'equivariante \`a gauche $\Phi ^{-1}(
\xi )\in C^{\infty}(G,TG)$, d\'efinie par $\Phi ^{-1}(\xi )(\gamma )=(T_{s(
\gamma )}L_{\gamma})\xi (s(\gamma ))$ et le crochet de Lie sur $C^{\infty}
(G^{(0)},{\mathcal G})$ est donn\'e par $\left[ \xi ,\eta \right] =\Phi \left(
\left[ \Phi ^{-1}(\xi ),\Phi ^{-1}(\eta )\right] \right)$. 

\medskip

On note $\Omega (\xi )=(\psi ')^{-1}\circ \Phi ^{-1}(\xi )\circ \psi$
la forme de $\Phi ^{-1}(\xi )$ dans les cartes.

\bigskip

\centerline{
\setlength \unitlength {0.25em}
\begin{picture}(150,30)(-75,-15)
\put (-30, -6){\vector(0,1){12}}
\put (40, -6){\vector(0,1){12}}
\put (-25, 10){\vector(1,0){59}}
\put (-21, -10){\vector(1,0){51}}
\put(-30, 10){\object{G}}
\put(-40, -10){\object{\R ^n\times \R ^m\supset U\times V}}
\put(40, 10){\object{TG}}
\put(40, -10){\object{T(U\times V)}}
\put(44, 0){\object{\psi '}}
\put(-34, 0){\object{\psi }}
\put(0, 13){\object{\Phi ^{-1}(\xi )}}
\put(0, -7){\object{\Omega (\xi )}}
\end{picture}}

\bigskip

On montre dans cette \'etape que 

\begin{equation}\label{Oomega}
\Omega (\xi )(u,v)=\left( 0,\xi _0(\sigma (u,v))+B(u,v,\xi _0(\sigma (u,v)))
\right)
\end{equation}

Par d\'efinition 

$\Phi ^{-1}(\xi )(\psi (u,v))=\left( T_{s(\psi (u,v))}L_{\psi 
(u,v)}\right) \xi (s(\psi (u,v)))$

\hspace*{1cm}$=\left( T_{s(\psi (u,v))}L_{\psi 
(u,v)}\right) \left( \dfrac{\d \psi}{\d v}(\sigma (u,v),0)\xi _0(\sigma (u,v)) 
\right)$

\hspace*{1cm}$=T_0\left( L_{\psi (u,v)}\psi (\sigma (u,v),\cdot )\right) \xi _0
(\sigma (u,v))$

\medskip

\noindent Mais $L_{\psi (u,v)}\psi (\sigma (u,v),w)=\psi (u,v)\psi (\sigma 
(u,v),w)=\psi (u,p(u,v,w))$, ce qui par d\'erivation conduit \`a 
$T_0\left( L_{\psi (u,v)}\psi (\sigma (u,v),\cdot )\right) =
\dfrac{\d \psi}{\d v}(u,v) \dfrac{\d p}{\d w}(u,v,0)$, donc

\begin{equation}\label{phiminus}
\Phi ^{-1}(\xi )(\psi (u,v))=\dfrac{\d \psi}{\d v}(u,v) \dfrac{\d p}{\d w}
(u,v,0)\xi _0(\sigma (u,v))
\end{equation}

\noindent Comme $p(u,v,w)=v+w+B(u,v,w)+O_3(u,v,w)$ on a 

\begin{equation}\label{derivp}
\dfrac{\d p}{\d w}(u,v,0)=I+\left(
\begin{array}{cccc}
B_1(u,v,f_1) & B_1(u,v,f_2) & ... & B_1(u,v,f_m)\\
\vdots & \vdots & \ddots & \vdots \\
B_m(u,v,f_1) & B_m(u,v,f_2) & ... & B_m(u,v,f_m)
\end{array}\right)
\end{equation}

\noindent En rempla\c cant \ref{derivp} dans \ref{phiminus} on obtient

\begin{equation}\label{phifinal}
\Phi ^{-1}(\xi )(\psi (u,v))=\dfrac{\d \psi}{\d v}(u,v)\xi _0(\sigma (u,v))+
\dfrac{\d \psi}{\d v}(u,v)B(u,v,\xi _0(\sigma (u,v)))
\end{equation}

\noindent La formule \ref{Oomega} r\'esulte comme suit :

\smallskip

$\Omega (\xi )(u,v)=\psi '(u,v)^{-1}\left( \Phi ^{-1}(\xi )(\psi (u,v))\right)$

\hspace*{1,9cm}$=\psi '(u,v)^{-1}\dfrac{\d \psi}{\d v}(u,v)\left( \xi _0(
\sigma (u,v))+B(u,v,\xi _0(\sigma (u,v)))\right)$

\hspace*{1,9cm}$=\left( 0,\xi _0(\sigma (u,v))+B(u,v,\xi _0(\sigma (u,v)))
\right)$

\medskip

{\bf 3.} La formule \ref{Oomega} montre en particulier que $\Omega 
(e_i )(u,v)=(0,f_i+B(u,v,f_i))$. Pour les \'el\'ements de cette forme le 
crochet de Lie dans $T(U\times V)$ est donn\'e par 

\medskip

$\left[ (0,a(u,v)),(0,b(u,v))\right] =\left( 0,\stackunder{k}{\dsum}\left( a_k(
u,v)\dfrac{\d b}{\d v_k}(u,v)-b_k(u,v)\dfrac{\d a}{\d v_k}(u,v)\right) \right)$

\medskip

\noindent En utilisant la bilin\'earit\'e de $B$ on obtient 
$\left[ \Omega (e_i),\Omega (e_j)\right] (u,v) =$

\noindent $=\left( 0,\stackunder{k}{\dsum}(\delta _{ik}+B_k(u,v,f_i))B(u,f_k,
f_j)-\stackunder{k}{\dsum} (\delta _{jk}+B_k(u,v,f_j))B(u,f_k,f_i) \right)$

\noindent $=\left( 0,B(u,f_i,f_j)+B(u,B(u,v,f_i),f_j)-B(u,f_j,f_i)-B(u,B(u,v,
f_j),f_i)\right)$

\medskip

\noindent et en particulier $\left[ \Omega (e_i),\Omega (e_j)\right] (u,0)=$

\noindent $=\left( 0,B(u,f_i,f_j)-B(u,f_j,f_i)+B(u,B(u,0,f_i),f_j)-
B(u,B(u,0,f_j),f_i)\right)$

\smallskip

\noindent $=\left( 0,B(u,f_i,f_j)-B(u,f_j,f_i)\right).$ 

\noindent On peut ainsi conclure

\medskip

$\left[ e_i,e_j\right] (\varphi (u))=\left[ \Phi ^{-1}(e_i),\Phi ^{-1}(e_j)
\right] (\psi (u,0))=\psi '(u,0)\left[ \Omega (e_i),\Omega (e_j)\right] (u,0)$

\smallskip

\hspace*{2,42cm}$=\dfrac{\d \psi}{\d v}(u,0)
\stackunder{k}{\dsum}\left( B_k(u,f_i,f_j)-B_k(u,f_j,f_i)\right) f_k$

\smallskip

\hspace*{2,42cm}$=\stackunder{k}{\dsum}\left( B_k(u,f_i,f_j)-B_k(u,f_j,f_i)
\right) e_k(\varphi (u))$. 

\hfill\QED

\begin{cor}\label{crochetLIE}
Pour toutes sections $\xi ,\eta \in C^{\infty}(G^{(0)},{\mathcal G})$ on a 

\smallskip

\centerline{$\left[ \xi ,\eta \right] _0(u)=B(u,\xi _0(u),\eta _0(u)-B(u,
\eta _0(u),\xi _0(u))$}
\end{cor}

{\it Preuve.}
On d\'eveloppe $\xi =\dsum \xi _ie_i$ et $\eta =\dsum \eta _je_j$ dans le 
rep\`ere mobile $\left\{ e_1,...,e_m\right\}$ et on s'en sert de la 
proposition pr\'ecedente pour remplacer $c_{ijk}$ dans $\left[ \xi ,\eta 
\right] =\dsum \xi _i\eta _jc_{ijk}e_k$. Il ne reste qu'\`a utiliser 
la bilin\'earit\'e de $B$. \QED

{\it Exemple.}
Soit $G$ un groupe de Lie d'unit\'e $e$. Dans ce cas $U=\left\{ 0\right\}$ et
$\psi :V\rightarrow G$ est une carte v\'erifiant $\psi (0)=e$. On associe \`a
$\psi$ la carte $\theta=\psi '(0):\R ^m\rightarrow {\mathcal G}$ de l'alg\`ebre de 
Lie. 
Par la proposition \ref{Taylorgroupoide} le produit dans $G$ est 
de la forme $\psi (v)\psi (w)=\psi (p(v,w))$,
o\`u $p:V\times V\rightarrow V$ est une application diff\'erentiable
qui admet un d\'eveloppement de la forme  $p(v,w)=v+w+B(v,w)+O_3(v,w)$ avec 
$B$ bilin\'eaire et l'inversion est donn\'ee par 
$\psi (v)^{-1}=\psi (-v+B(v,v)+O_3(v))$.
Comme $\sigma=0$ on retrouve $a_{ij}=0$ et la proposition \ref{329} montre que
les constantes de structure $c_{ijk}$ de l'alg\`ebre de Lie sont donn\'ees par 
$c_{ijk}=B_k(f_i,f_j)-B_k(f_j,f_i)$. 

On retrouve ainsi des r\'esultats connus pour les groupes et alg\`ebres de 
Lie qu'on peut trouver par exemple dans {\em \cite{kirillov}}.

\smallskip

{\it Exemple.}
Soit $G=M\times M$, le groupo{\"\i}de principal transitif associ\'e \`a la 
vari\'et\'e $M$, un \'el\'ement $(x,x)\in G^{(0)}$ et $\alpha :U_1\rightarrow 
M$ une carte de $M$ telle que $\alpha (0)=x$. On peut prendre la carte 
$\psi :U\times V\rightarrow G$, donn\'ee par $\psi (u,v)=(\alpha (u),\alpha 
(u+v))$, o\`u $U,V$ sont des voisinages de $0\in \R ^m$ telles que $U\in U_1$ 
et $U+V\in U_1$. La carte de $TM$ engendr\'ee par $\psi$ est $\theta :U\times 
\R ^m\rightarrow {\mathcal G}$, $\theta (u,v)=(\alpha (u),\alpha '(u)v)$. Dans ce 
cas $\sigma (u,v)=u+v$, $p(u,v,w)=v+w$, $B(u,v,w)=0$. Pour les fonctions
de structure de $TM$ on retrouve $a_{ij}=\delta _{ij}$ et $c_{ijk}=0$.

\noindent {\it Remerciements} Les r\'esultats de cet article ont 
\'et\'e obtenu durant mon sejour \`a l'Universit\'e d'Orl\'eans et 
je tiens a exprimer ma gratitude \`a Jean Renault
pour tout son soutien et ses tr\`es judicieuses remarques
et \`a Claire Anantharaman qui  
au d\'ebut de ce travail m'a fait comprendre le cas des groupes de Lie.

{\small
\bibliographystyle{french}

}

\bigskip
{\hspace*{8.5cm}Acad\'emie Roumaine

\hspace*{8cm}Institut de Math\'ematiques 

\hspace*{7.5cm}Calea Grivi\c tei 21, P.O. Box 1-764, 

\hspace*{8cm}Bucarest 70700, Romania}

\hspace*{6.6cm}{E-mail address}: ramazan@pompeiu.imar.ro

\end{document}